\documentclass{amsart}

\usepackage{amsfonts}
\usepackage{amsmath}
\usepackage{amsthm}
\usepackage{amssymb}
\usepackage[english]{babel}
\usepackage{graphics}
\usepackage{latexsym}
\usepackage{longtable} 
\usepackage{mathrsfs} 
\usepackage{graphicx}
\usepackage{wrapfig} 
\usepackage[pdftex,colorlinks=true,linkcolor=blue,urlcolor=blue,citecolor=blue]{hyperref}

\newtheorem{theorem}{Theorem}
\newtheorem{corollary}[theorem]{Corollary}
\newtheorem{lemma}[theorem]{Lemma}

\newtheorem{claim}[theorem]{Claim}
\newtheorem{example}[theorem]{Example}

\theoremstyle{definition}
\newtheorem{definition}[theorem]{Definition}
\newtheorem{remark}[theorem]{Remark}

\newcommand{\mA}{\mathcal{A}}

\newcommand{\mP}{\mathscr{P}}
\newcommand{\mM}{\mathcal{M}}

\newcommand{\mY}{\mathscr{Y}}
\newcommand{\mD}{\mathcal{D}}

\newcommand{\A}{\textbf{A}}
\newcommand{\R}{\mathbb{R}}

\newcommand{\N}{\mathbb{N}}

\newcommand{\mB}{\mathbb{B}}
\newcommand{\X}{\textbf{X}}
\newcommand{\Y}{\textbf{Y}}

\newcommand{\noi}{\noindent}
\newcommand{\ms}{\medskip}

\newcommand{\al}{\alpha}
\newcommand{\be}{\beta}
\newcommand{\ga}{\gamma}

\newcommand{\de}{\delta}
\newcommand{\De}{\Delta}
\newcommand{\e}{\varepsilon}

\newcommand{\la}{\lambda}

\newcommand{\Om}{\Omega}

\newcommand{\lharpoonup}{-\!\!\!\!-\!\!\!\!\rightharpoonup}

\newcommand{\weak }{\, -\!\!\!\!-\!\!\!\rightharpoonup}
\newcommand{\weakstar }{ \overset{\, *_{\phantom{|}}}{{\smash{\weak }}\, } }

\newcommand{\larrow}{\longrightarrow}
\newcommand{\ot}{\otimes}

\newcommand{\lmapsto}{\longmapsto}
\newcommand{\ri}{\rightarrow}

\newcommand{\p}{\partial}
\newcommand{\sub}{\subseteq}
\newcommand{\set}{\setminus}
\newcommand{\by}{\times}

\newcommand{\ess}{\textrm{ess}}

\newcommand{\dist}{\textrm{dist}}

\newcommand{\supp}{\textrm{supp}}

\newcommand{\bt}{\begin{theorem}}\newcommand{\et}{\end{theorem}}
\newcommand{\bd}{\begin{definition}}\newcommand{\ed}{\end{definition}}
\newcommand{\bl}{\begin{lemma}}\newcommand{\el}{\end{lemma}}
\newcommand{\beq}{\begin{equation}}\newcommand{\eeq}{\end{equation}}
\newcommand{\bc}{\begin{claim}}\newcommand{\ec}{\end{claim}}
\newcommand{\bex}{\begin{example}}\newcommand{\eex}{\end{example}}
\newcommand{\bcor}{\begin{corollary}}\newcommand{\ecor}{\end{corollary}}
\newcommand{\bp}{\begin{proof}}\newcommand{\ep}{\end{proof}}

\newcommand{\BPL}{\medskip \noindent \textbf{Proof of Lemma} }

\newcommand{\BPCOR}{\medskip \noindent \textbf{Proof of Corollary} }

\newcommand{\BPT}{\medskip \noindent \textbf{Proof of Theorem} }

\numberwithin{equation}{section}

\begin{document}

\title[A New Variational Characterisation of $\infty$-Harmonic Maps]{A New Characterisation of $\infty$-Harmonic and $p$-Harmonic Maps via Affine Variations in $L^\infty$}

\author{Nikos Katzourakis}

\thanks{\texttt{The author has been financially supported by the EPSRC grant EP/N017412/1.} }

\address{Department of Mathematics and Statistics, University of Reading, Whiteknights, PO Box 220, Reading RG6 6AX, Berkshire, UK}
\email{n.katzourakis@reading.ac.uk}



\date{}



\begin{abstract} Let $u: \Omega \subseteq \mathbb{R}^n \longrightarrow \mathbb{R}^N$ be a smooth map and $n,N \in \mathbb{N}$. The $\infty$-Laplacian is the PDE system
\[ \tag{1} \label{1}
\Delta_\infty u \, :=\, \Big(Du \otimes Du + |Du|^2[Du]^\bot\!  \otimes I\Big) :D^2u\, =\, 0,
\]
where $[Du]^\bot := \text{Proj}_{R(Du)^\bot}$. \eqref{1} constitutes the fundamental equation of vectorial Calculus of Variations in $L^\infty$, associated to the model functional
\[ \tag{2} \label{2}
E_\infty (u,\Omega')\, =\, \big\| |Du|^2\big\|_{L^\infty(\Omega')} ,\ \ \ \Omega' \Subset \Omega.
\]
We show that generalised solutions to \eqref{1} can be characterised in terms of \eqref{2} via a set of designated affine variations. For the scalar case $N=1$, we utilise the theory of viscosity solutions of Crandall-Ishii-Lions. For the vectorial case $N\geq 2$, we utilise the recently proposed by the author theory of $\mathcal{D}$-solutions. Moreover, we extend the result described above to the $p$-Laplacian, $1<p<\infty$. 
\end{abstract}

\maketitle


\section{Introduction} \label{section1}

Let $n,N\in \N$. Given a (smooth) map $u:\Om \sub \R^n \larrow \R^N$ defined on an open set, let $\R^{Nn}$ and $\R^{Nn^2}_s$ denote respectively the space of matrices and the space of symmetric tensors wherein the gradient matrix and the hessian tensor
\[
Du(x)\ =\ \big( D_i u_\al(x)\big)_{i=1,...,n}^{\al=1,...,N},\ \ \ \ D^2 u(x)\ =\ \left( D^2_{ij}u_\al(x)\right)_{i,j =1,...,n}^{\al=1,...,N}
\]
of  $u$ are valued. Obviously, $D_i\equiv \p/\p x_i$, $x=(x_1,...,x_n)^\top$, $u=(u_1,...,u_N)^\top$. 
In this paper we are primarily interested in the so-called $\infty$-Laplacian which is the following quasilinear $2$nd order nondivergence system:
\beq \label{1.1}
\De_\infty u \, :=\, \Big(Du \ot Du + |Du|^2[Du]^\bot \! \ot I \Big):D^2u\, =\, 0.
\eeq
Here $[Du]^\bot$ denotes the orthogonal projection on the orthogonal complement of the range of $Du$ and $|Du|$ is the Euclidean norm of $Du$ on $\R^{Nn}$. In index form \eqref{1.1} reads  
\[
\begin{split}
\sum_{\be=1}^N\sum_{i,j=1}^n \Big(D_i u_\al \, D_ju_\be + & \, |Du|^2  [Du]_{\al \be}^\bot \, \de_{ij}\Big)\, D_{ij}^2u_\be\, =\, 0, \ \ \ \ \al=1,...,N, \\
& [Du]^\bot\, :=\, \text{Proj}_{(R(Du))^\bot}.
\end{split}
\]
We are also interested in the more classical $p$-Laplacian for $1<p<\infty$, which is the following divergence system:
\beq \label{1.2}
\De_p u \, :=\, \text{div}\big(|Du|^{p-2}Du \big)\, =\, 0.
\eeq
The system \eqref{1.1} is the fundamental equation which arises in vectorial Calculus of Variations in the space $L^\infty$, that is in connection to variational problems for the model functional
\beq   \label{1.3}
E_\infty (u,\Omega')\, :=\, \big\| |Du|^2 \big\|_{L^\infty(\Omega')} ,\ \ \ \Omega' \Subset \Omega, \ \  u \in W^{1,\infty}_{\text{loc}}(\Om,\R^N).
\eeq
The scalar counterpart of \eqref{1.1} when $N=1$ simplifies to
\[
Du \ot Du :D^2u\, = \sum_{i,j=1}^n D_i u \, D_ju \, D^2_{ij}u\,=\, 0
\]
and first arose in the work of G.\ Aronsson in the 1960s (\cite{A1, A2} and for a pedagogical introduction see \cite{C,K7}) who pioneered the field of Calculus of Variations in the space $L^\infty$. The full system \eqref{1.1} first appeared in recent work of the author \cite{K1} who initiated the systematic study of the vectorial case in a series of papers \cite{K1}-\cite{K6} (see also the recent joint contributions with Abugirda, Ayanbayev, Croce, Pisante, Manfredi, Moser and Pryer \cite{AK, AyK, CKP, KP, KM, KM2, KP2}). Let us note also the early vectorial contributions by Barron-Jensen-Wang \cite{BJW1, BJW2} who, among other deep results, proved existence of absolute minimisers for general supremal functionals in the ``rank-one" cases $\min\{n,N\}=1$ and also defined and studied the correct vectorial $L^\infty$-version of quasiconvexity. However, their fundamental contributions were at the level of the functional and the correct (non-obvious) vectorial counterpart of Aronsson's equation was not known at the time.

On the other hand, the $p$-Laplacian \eqref{1.2} is a classical model which arises in conventional Calculus of Variations for integral functionals, in particular as the Euler-Lagrange equation of 
\beq   \label{1.4}
E_p (u,\Omega')\, :=\, \big\| |Du|^p \big\|_{L^1(\Omega')} ,\ \ \ \Omega' \Subset \Omega, \ \  u \in W^{1,p}_{\text{loc}}(\Om,\R^N).
\eeq

A standard difficulty in both the scalar and the vectorial case of \eqref{1.1} is that it is nondivergence and since in general smooth solutions do not exist, the definition of generalised solutions is an issue. In the vectorial case, an additional difficulty is that the system has discontinuous coefficients even if the solution might be smooth  (see \cite{K2}). This happens because the projection $[Du(x)]^\bot$ ``feels" the dimension of the tangent space $R(Du(x))\sub \R^N$.

In this paper we are concerned with the variational characterisation of appropriately defined generalised solutions to \eqref{1.1} and \eqref{1.2} in both the scalar and the vectorial case in terms of the supremal functional \eqref{1.3}. The main results of this paper are contained in the statements of Theorems \ref{theorem8}, \ref{theorem10} and \ref{theorem11} (and Corollaries \ref{corollary9}, \ref{corollary12}). Roughly speaking, these results claim that for $1<p\leq \infty$  we have
\[
\De_p u\, =\, 0 \text{ on }\Om \ \ \Longleftrightarrow\ \ 
\left\{
\begin{array}{l}
\text{For all }\, \Om'\Subset \Om \text{ and }  A\in \mA^p_{\Om'}(u), \ms\\
 E_\infty(u,\Om') \,\leq\, E_\infty(u+A,\Om')
 \end{array}
 \right.
\]
where $\mA^p_{\Om'}(u)$ is a designated set of \textbf{affine} mappings depending on $u$ and on the subdomain $\Om'$. This result is quite surprising in that both the $\infty$-Laplacian \eqref{1.1} and the $p$-Laplacian \eqref{1.2}  are associated to the respective supremal/integral functionals \eqref{1.3},  \eqref{1.4} (and not both associated to \eqref{1.3}) when the classes of variations are \textit{compactly supported}. In the scalar case, the appropriate notion of minimisers for \eqref{1.3} characterising $\infty$-Harmonic functions has been discovered by Aronsson and today we know several more characterisations involving e.g.\ comparison, Lipschitz extensions and Game Theory (see \cite{C,K7}). In the vectorial case, the correct extension of Aronsson's notion of Absolute Minimals which characterises \eqref{1.1} via \eqref{1.3} has been introduced in \cite{K4}.

A central point in both the statements and the proofs of our main results Theorems \ref{theorem8}, \ref{theorem10} and \ref{theorem11} is that solutions to \eqref{1.1}-\eqref{1.2} in general are nonsmooth and they need to be considered in a generalised sense. We  discuss below about generalised solutions separately when $N=1$ and $N\geq 2$.

For the scalar case, we invoke the well established notion of viscosity solutions of Crandall-Ishii-Lions \cite{CIL} which effectively is based on the maximum principle. Since the $p$-Laplacian is singular for $1<p<2$, we actually use a ``feeble" variant of the original viscosity notions taken from \cite{K0}. Although \eqref{1.2} is in divergence from and the natural definition of weak solution to it is via duality, we find it more fruitful to treat it instead in the viscosity sense. Due to the results in the aforementioned papers, it is known that viscosity and weak solutions of the $p$-Laplacian coincide.

For the vectorial case, things are much more intricate. Motivated by \eqref{1.1}, in the very recent works \cite{K9, K8} we introduced a new duality-free theory of weak solutions which allows for just measurable maps to be rigorously interpreted and studied as solutions to PDE systems of any order 
\beq \label{1.5}
F\Big(\cdot,u,Du,D^2u,...,D^pu\Big)\, =\, 0\ \ \ \text{ on } \Om,
\eeq
which can be allowed to have even discontinuous coefficients. Using this new approach, in the  papers \cite{K8}-\cite{K11} we studied efficiently certain problems which we discuss briefly at the end of the introduction.

Our generalised solutions are not based either on integration-by-parts or on the maximum principle. Instead, we build on the probabilistic interpretation of limits of difference quotients by utilising \emph{Young} measures valued into compactifications. We caution the reader that we are not using the ``standard" Young measures of Calculus of Variations and of PDE theory which are valued into Euclidean spaces (see e.g.\ \cite{E, P, FL, CFV, FG, V, KR}). In the current setting, Young measures valued into spheres are utilised by applying them to the difference quotients of our candidate solution. The motivation for $W^{1,\infty}_{\text{loc}}$ solutions of $2$nd order systems which are relevant to this paper is the following: let $u\in W^{2,\infty}_{\text{loc}}(\Om,\R^N)$ be a strong solution to a $2$nd order system of the form
\beq   \label{1.6}
F\big(Du(x),D^2u(x) \big)\,=\, 0, \quad \text{ a.e.\ }x\in \Om.
\eeq
We now rewrite \eqref{1.6} in the following unconventional fashion
\[
\sup_{\X_x\in \, \supp(\de_{D^2 u(x)})} \big| F\big(Du(x),\X_x\big)\big|\, =\, 0, \quad \text{ a.e. }x\in \Om
\]
and we view the hessian $D^2u $ as a probability-valued mapping given by the Dirac mass: $\de_{D^2u}$. The hope is then that we may relax the requirement to have concentration measures and allow instead general probability-valued maps arising as limits of difference quotients of $W^{1,\infty}_{\text{loc}}$ maps. Indeed, if $u : \Om \sub \R^n \larrow \R^N$ is just $W^{1,\infty}_{\text{loc}}$, we may view the usual difference quotients of $Du$ as Young measures into the 1-point compactification
\[
\de_{D^{1,h}Du} \ : \ \Om \sub \R^n \larrow \mathscr{P}\big( \smash{\overline{\R}}^{Nn^2}_s \big) , \ \ \ x\lmapsto \de_{D^{1,h}Du(x)}
\]
(see Section \ref{section2} for the precise definitions). Since the Young measures are a weakly* compact set, there exist probability-valued limit maps such that along infinitesimal subsequences $(h_\nu)_1^\infty$ we have
\beq  \label{1.7}
\de_{D^{1,h_\nu}Du} \, \weakstar \, \mD^2 u, \ \ \text{ in the Young measures, as }\nu \ri \infty
\eeq
(even if $u$ is merely $W^{1,\infty}_{\text{loc}}$). Then, we require
\beq  \label{1.8}
\sup_{\X_x\in \supp(\mD^2 u(x))\set \{\infty\}} F\big(Du(x),\X_x\big)\, =\, 0, \quad  \ \text{ a.e. }x\in \Om,
\eeq
 for any``diffuse hessian" $\mD^2 u$.  Since \eqref{1.7} and \eqref{1.8} are \emph{independent} of the twice differentiability of $u$, they can be taken as a notion of generalised solution which we call $\mD$-solutions. In the event that $u \in W^{2,\infty}_{\text{loc}}$, then $\mD^2 u=\de_{D^2u}$ and we reduce to strong solutions. 
  
A flaw of our characterisations is that we require our generalised solutions to be $C^1$ and not just $W^{1,\infty}_{\text{loc}}$. This is not a restriction for the $p$-Laplacian since it is well know that $p$-Harmonic maps are $C^{1,\al}$ (\cite{U}). However, except for the case of $n=2$, $N=1$ (see Savin and Evans-Savin \cite{S,ES}), the $C^1$ regularity of $\infty$-Harmonic functions (and a fortiori of maps) is an open problem, at least to date. However, even with the extra $C^1$ hypothesis, the results are new even in the scalar case. We believe that they are interesting anyway and might allow to glean more information that will unravel the still largely mysterious behaviour of $\infty$-Harmonic functions (and maps). For the $p$-Laplacian we restrict our attention only to $N=1$ and we refrain from extending Theorem \ref{theorem10} to $N\geq2$. This however can be done relatively easily along the lines of Theorem \ref{theorem11}. 

Further, we postpone the discussion of the more difficult question of relation of viscosity and $\mD$-solutions for future work. It is easily seen though that $\mD$-solutions do not have comparison built in the notion as viscosity solutions (in the vectorial case in general not even $C^\infty$-solutions are unique, see \cite{K5}) and hence $\mD$-solutions are not stronger than viscosity solutions. On the other hand, absolutely minimising $\mD$-solutions are viscosity solutions and we conjecture that the opposite is true as well. (Let us note that in \cite{KP2} is was recently proved that absolutely minimising $\mD$-solutions of higher order $L^\infty$ variational problems are unique.)

We conclude this introduction with certain interesting results we have obtained via the new theory of $\mD$-solutions. In the paper \cite{K8} we proved existence to the Dirichlet problem for \eqref{1.1} (uniqueness of smooth solutions has been disproved in \cite{K5}). Again in \cite{K8}, we also proved uniqueness and existence to the Dirichlet problem for the fully nonlinear degenerate elliptic system $F(\cdot,D^2u)=f$. In \cite{K9} we proved existence to the Dirichlet problem for the system arising from the functional
\[
\ \ \ I_\infty(u,\Om')\, :=\, \big\| H(\cdot,u,u') \big\|_{L^\infty(\Om')}, \ \ \ \ \ u \,:\ \Om\sub \R\larrow \R^N,\ \Om'\Subset\Om.
\] 
In \cite{K10} we established the equivalence between weak and $\mD$-solutions to linear symmetric hyperbolic systems and in \cite{K11} we developed a systematic mollification method for $\mD$-solutions. We finally note that to the best of our knowledge, the only vectorial contribution by other authors relevant to the content of this paper is the work of Sheffield-Smart \cite{SS} which however is restricted to the class of smooth solutions.

\section{Basics on generalised solutions to fully nonlinear systems} \label{section2}

We begin with some basic material. A much more detailed introduction of the theory of $\mD$-solutions can be found in \cite{K8}-\cite{K11}.

\ms

\noi \textbf{Preliminaries.} Let  $u:\Om\sub \R^n \larrow \R^N$ be a map defined over an open set. Unless indicated otherwise, Greek indices $\al,\be,\ga,...$ will run in $\{1,...,N\}$ and latin indices $i,j,k,...$ will run in $\{1,...,n\}$. The norm symbols $|\cdot|$ will always mean the Euclidean ones, whilst Euclidean inner products will be denoted by either ``$\cdot$" on $\R^n,\R^N$ or by ``$:$" on $\R^{Nn},\R^{Nn^2}_s$. For example, 
\[
|\X|^2 \, =\, \X :\X \, =\, \sum_{\al=1}^N \sum_{i, j=1}^n \X_{\al i j} \X_{\al i j}, \ \ \ \ \X\in \R^{Nn^2}_s,
\]
etc. Our measure theoretic and function space notation is either standard as e.g.\ in \cite{E,E2} or self-explanatory. For example, ``measurable" means ``Lebesgue measurable", the Lebesgue measure will be denoted by $|\cdot|$, the $L^p$ spaces of maps $u$ as above by $L^p(\Om,\R^N)$, etc. Especially for the space $L^\infty(\Om,\R^{Nn})$, we will simplify the notation and since the norm on $\R^{Nn}$ is always the Euclidean, we will write 
\[
\|Du\|_{L^\infty(\Om)}\, =\, \underset{\Om}{\ess\,\sup} \, |Du| .
\]
We will systematically use the Alexandroff $1$-point compactification of the space $\R^{Nn^2}_s$. Its topology will be the one which makes it homeomorphic to the sphere of dimension $Nn(n+1)/2$ (via the stereographic projection which identifies the north pole with $\{\infty\}$). We will denote it by 
\[
\smash{\overline{\R}}^{Nn^2}_s\, :=\ \R^{Nn^2}_s \cup \{\infty\}.
\]
Then, the space $\R^{Nn^2}_s$ will be viewed as a metric vector space, isometrically contained into its $1$-point compactification.

\ms

\noi \textbf{Young Measures.} Let $\Om\sub \R^n$ be open. The Young measures can be identified with a subset of the unit sphere of a certain $L^\infty$ space of measure-valued maps and this provides very useful properties, like compactness.

\begin{definition} The set of Young Measures from $\Om$ to $\smash{\overline{\R}}^{Nn^2}_s$ is the subset of the unit sphere of the space $L^\infty_{w^*}\big( \Om,\mM\big(\smash{\overline{\R}}^{Nn^2}_s\big) \big)$ which contains probability-valued maps:
\[
\mY\big(\Om,\smash{\overline{\R}}^{Nn^2}_s\big)\, :=\, \Big\{ \vartheta\, \in \, L^\infty_{w^*}\big( \Om,\mM\big(\smash{\overline{\R}}^{Nn^2}_s\big) \big)\, : \, \vartheta(x) \in \mP \big(\smash{\overline{\R}}^{Nn^2}_s\big),\text{ for a.e. }x\in \Om\Big\}.
\]
\end{definition}
The space $L^\infty_{w^*}\big( \Om, \mM\big(\smash{\overline{\R}}^{Nn^2}_s\big) \big)$ is a dual Banach space and consists of measure-valued maps $\Om \ni x \longmapsto \vartheta(x)  \in \mM\big(\smash{\overline{\R}}^{Nn^2}_s\big)$ which are weakly* measurable, in the sense that for any Borel set $ \mathcal{U} \sub \smash{\overline{\R}}^{Nn^2}_s $, the function $\Om\ni x \lmapsto [\vartheta(x)]( \mathcal{U} )  \in \R$ is measurable. The norm of the space is given by
\[
\| \vartheta \|_{ L^\infty_{w^*} ( \Om, \mM (\smash{\overline{\R}}^{Nn^2}_s) ) }\, :=\, \underset{x\in \Om}{\ess\,\sup}\, \left\|\vartheta(x) \right\| \big( \smash{\overline{\R}}^{Nn^2}_s \big)
\]
where ``$\|\cdot\|$" denotes the total variation. For background material on these spaces we refer e.g.\ to \cite{FL, Ed, V} and to \cite{K8}-\cite{K11}. The $L^\infty_{w^*}$ space above is the dual space of the space $L^1 \big( \Om, C^0\big(\smash{\overline{\R}}^{Nn^2}_s\big) \big)$ of Bochner integrable maps. The points of this $L^1$ space are the Carath\'eodory functions $\Phi : \Om \by \smash{\overline{\R}}^{Nn^2}_s \larrow \R$ which satisfy
\[
\| \Phi \|_{L^1 ( \Om, C^0 (\smash{\overline{\R}}^{Nn^2}_s))}\, :=\, \int_\Om  \big\| \Phi(x,\cdot)\big\|_{C^0(\R^{Nn^2}_s )} \, dx\, <\, \infty.
\]
It is well known that the unit ball of $L^\infty_{w^*}$ is sequentially weakly* compact. Hence, for any bounded sequence $(\vartheta^m)_1^\infty \sub L^\infty_{w^*} $, there is a limit map $\vartheta$ and a subsequence of $m$'s along which $\vartheta^m \overset{^*}{\smash{\lharpoonup}}\, \vartheta$ as $m\ri \infty$.

\begin{remark}[Properties of Y.M.] \label{remark2} We note the following facts about Young measures (proofs can be found e.g.\ in \cite{FG}): 

\smallskip
\noi i)  [\textit{Functions as Y.M.}]  The set of measurable maps $U : \Om\sub \R^n \larrow \smash{ {\R}}^{Nn^2}_s$ can be identified with a subset of the Young measures via the embedding $U \mapsto \de_U$, $\de_U(x):= \de_{U(x)}$. 

\smallskip
 
\noi ii) [\textit{Weak* compactness of Y.M.}] The set of Young measures is convex and sequentially compact in the weak* topology induced from $L^\infty_{w^*}$.  

\end{remark}

The next lemma is a minor variant of a classical result (see \cite{FG, FL,K8}) but it plays a fundamental role in our setting because it guarantees the compatibility of classical/strong solutions with $\mD$-solutions. 

\begin{lemma} \label{lemma0} Let $U^\nu,U^\infty : \Om\sub \R^n\larrow \smash{{\R}}^{Nn^2}_s$ be measurable maps, $\nu\in \N$. Then, up the passage to a subsequence, the following equivalence holds true:
\[
 \de_{U^\nu}  \weakstar \de_{U^\infty} \text{ in }\mY\big(\Om, \smash{\overline{\R}}^{Nn^2}_s\big)  \ \ \ \Longleftrightarrow \ \ \ U^{\nu} \larrow U^\infty \ \text{ a.e.\ on }\Om.\ \ \ \
 \]
\end{lemma}

\smallskip

\noi \textbf{The notion of $\mD$-Solutions to fully nonlinear $2$nd order systems.}  Herein we consider the special case of once differentiable solutions to second order systems which is relevant to the $\infty$-Laplacian. For the general case of measurable solutions to $p$th order system we refer to \cite{K8,K11}.

Let $D^{1,h}$ denote the usual difference quotient operator on $\R^n$, i.e.\ given a map $v : \Om \sub \R^n\larrow \R^N$ and $h\neq 0$, we understand $v$ as being extended by zero on $\R^n\set\Om$ and we set
\[
\begin{split}
D^{1,h}_i v(x)\, &:=\, \frac{v(x+he^i) - v(x)}{h} , \ \ \ \ \ \ \ \ \ \ \ x\in \Om, \\ 
 D^{1,h}v(x)\, &:=\, \left(D^{1,h}_1v(x),...,D^{1,h}_n v (x)\right), \ \ x\in \Om.
\end{split}
\]
\begin{definition}[Diffuse Hessians]  \label{Diffuse hessians} Let $u : \Om\sub \R^n\larrow\R^N$ be a locally Lipschitz continuous map. We define the \textbf{diffuse hessians $\mD^2 u$ of $u$} as the subsequential weak* limits of the difference quotients of the gradient in the space of Young measures along infinitesimal sequences $(h_{\nu})_{1}^\infty$:
\[
\de_{D^{1,h_{\nu_k}} Du} \weakstar \mD^2u \ \ \   \text{ in }\mY\big(\Om,  \smash{\overline{\R}}^{Nn^2}_s \big), \ \text{ as }k \ri \infty.
\]
\end{definition} 

Next is our notion of generalised solution for the vectorial case. We will use the notation ``$\supp_*$" to denote the \emph{reduced support} of a probability measure $\vartheta$ on $\smash{\overline{\R}}^{Nn^2}_s$ ``off infinity", namely,
\[ 
\supp_* (\vartheta ) \, := \, \supp  ( \vartheta ) \set \{\infty\}, \ \ \ \ \vartheta \in \mP\big( \smash{\overline{\R}}^{Nn^2}_s\big).
\]
\begin{definition}[Lipschitz $\mD$-solutions to $2$nd order systems] \label{definition13} Let $\Om\sub \R^n$ be an open set and
\[
F\ : \ \ \R^{Nn} \! \by \R^{Nn^2}_s  \larrow \R^N
\]
a mapping which is Borel measurable with respect to the first argument and continuous with respect to the second argument. Consider the PDE system
\beq \label{2.11a}
F\big(Du,D^2u\big)\, =\, 0\ \ \text{ on }\Om.
\eeq
We say that the locally Lipschitz continuous map $u : \Om\sub \R^n \larrow \R^N$ is a \textbf{$\mD$-solution of \eqref{2.11a}} when for any diffuse hessian $\mD^2 u$ of $u$, we have
\beq \label{2.11aa}
\sup_{ \X_x  \in \, \supp_* (\mD^2 u(x)) }\,
\left| F\big(Du(x), \X_x \big)\right| \, =\, 0, \quad \text{ a.e. }x\in\Om.
\eeq
\end{definition}
In particular, for the $\infty$-Laplace system \eqref{1.1}, a $W^{1,\infty}_{\text{loc}}$ map $u : \Om\sub \R^n \larrow \R^N$ is $\infty$-Harmonic in the  $\mD$-sense, when for a.e.\ $x\in\Om$ and all $\X_x  \in  \supp_* (\mD^2 u(x))$, we have
\[
\Big(Du(x) \ot Du(x) +|Du(x)|^2[Du(x)]^\bot \!\ot I \Big): \X_x \, =\, 0. 
\]
Note that at certain points it may happen that $\mD^2 u(x)=\de_{\{\infty\}}$ which implies that the reduced support of $\mD^2 u(x)$ is empty. The criterion then is understood to be trivially satisfied. Further, the $\mD$-notions are compatible with the strong/classical notions of solution: this is a direct consequence of Lemma \ref{lemma0} and the definition of diffuse hessians. 

\begin{remark}[An alternative formulation of $\mD$-solutions] \label{remark6} We give an alternative ``integral" form of Definition \ref{definition13} above which we put foremost in \cite{K8}-\cite{K10} because of its technical convenience for the existence/uniqueness proofs therein. We will not use this version herein, however. Note first that \eqref{2.11aa} can be rephrased as the following differential inclusion for the support:
\[
\ \ \supp(\mD^2 u(x)) \sub \Big\{ \X \in \R^{Nn^2}_s :\, \left| F\big(Du(x), \X \big)\right| =0 \Big\} \cup \{\infty \} ,\ \ \text{ a.e. }x\in \Om.
\]
Then, for any compactly supported $\Phi \in C^0_c\big( \R^{Nn^2}_s \big)$ off infinity and for a.e.\ $x\in \Om$, the continuous function
\[
\smash{\overline{\R}}^{Nn^2}_s \, \ni \X \lmapsto \Phi(\X)\,  F\big(Du(x), \X \big)  \in \R^N
\]
is well-defined on the compactification and also vanishes on the support of any diffuse hessian measure. As a consequence, we have the statement
\beq \label{2.12a}
\int_{  \smash{\overline{\R}}^{Nn^2}_s } \Phi(\X)\, F\big(Du(x), \X \big)\, d[\mD^2u(x)](\X)\, =\, 0, \ \ \ \text{ a.e. }x\in \Om,
\eeq
for any $\Phi \in C^0_c\big( \R^{Nn^2}_s \big)$ and any diffuse hessian $\mD^2u \in \mY\big(\Om,  \smash{\overline{\R}}^{Nn^2}_s \big)$. It can be easily seen that the converse is true as well (see \cite{K8}) and hence \eqref{2.12a} is a restatement of \eqref{2.11aa}.
\end{remark}

For more details on the material of this section (e.g.\ analytic properties, equivalent formulations of Definition \ref{definition13}, etc) we refer to \cite{K8}-\cite{K11}.

\ms

\noi \textbf{The notion of Feeble Viscosity Solutions to fully nonlinear $2$nd order equations.}  The definitions of this paragraph are taken from \cite{K0} (see also \cite{JJ,JLM} where the ``feeble" counterparts of the ``usual" viscosity notion first appeared) but here we apply them only to the case of the $p$-Laplacian for $1<p<\infty$. The standard viscosity notions as in  \cite{CIL,C,K7}  do not apply here because we treat also the singular case of the $p$-Laplacian when $p<2$ which is not even defined when the gradient vanishes.

Let $F : (\R^n\!\set\!\{0\}) \by\R^{n^2}_s \larrow \R$ be a continuous function which satisfies the monotonicity hypothesis $F(P,\X) \leq F(P,\Y)$ when $\X\leq \Y$ in $\R^{n^2}_s$. We consider the PDE
\[
F\big(Du,D^2u\big)\, = \, 0\ \ \text{ on }\Om.
\]  
Let $u : \Om\sub \R^n \larrow \R$ be a continuous function. Given a triplet $(x,P,\X)\in \Om \by \R^n \by\R^{n^2}_s$, we define the quadratic polynomial $T_{P,\X,x}u$  by setting
\[
T_{P,\X,x}u(z)\,  :=\, u(x)\, +\, P\cdot z \, +\, \frac{1}{2}\X:z\ot z,\ \ \ z\in \R^n.
\]
We then set
\[
J_0^{2,\pm} u(x) := \left\{(P,\X)\in\! (\R^n \! \set \! \{0\}) \by\R^{n^2}_s  \Bigg|  u(z+x) 
{
\! \begin{array}{c}  \leq \\ \geq \end{array} \!
}
T_{P,\X,x}u(z)+o(|z|^2), \text{ as }z\ri0 \right\}
\]
and call $J_0^{2,\pm} u(x)$ the \emph{feeble $2$nd order sub/superjet of $u$ at $x$}. We say that $u$ is a feeble viscosity solution of $F\big(Du,D^2u\big)\geq 0$ (resp.\ of $F\big(Du,D^2u\big)\leq 0$) on $\Om$ when for any $x\in \Om$
\[
\inf_{(P,\X) \in J_0^{2,+}u(x)} \, F(P,\X)\, \geq\, 0 \ \ \ \left( \text{ resp. }\sup_{(P,\X) \in J_0^{2,-}u(x)} \, F(P,\X)\, \leq\, 0 \right).
\]
Feeble viscosity solutions of $F\big(Du,D^2u\big)=0$ are defined as the combination of the above one-sided sub/super solution statements. 

If $u\in C^1(\Om)$, then any pair $(P,\X)$ in $J_0^{2,\pm}u(x)$ satisfies $P=Du(x)$. In this case we will use the notation
\[
D^{2,\pm}u(x)\, :=\, \Big\{ \X \in \R^{n^2}_s \ \Big| \ (Du(x),\X) \in J_0^{2,\pm}u(x)\Big\}
\]
and we will call $D^{2,\pm}u(x)$ the set of feeble $2$nd order sub/super derivatives of $u$ at $x\in \Om$.

\section{Two elementary lemmas}

In this brief section we isolate a couple of very simple technical results which contain an essential common part of the proofs of the main results in both the scalar and the vectorial case.

\begin{lemma} \label{lemma1} Let $\Om\sub \R^n$ be open and $u\in C^1(\Om,\R^N)$.  Given $\Om'\Subset \Om$, we set
\[
\Om'(u)\, :=\, \Big\{ x\in \overline{\Om'}\ \Big| \ |Du(x)|\, =\, \|Du\|_{L^\infty(\Om')}\Big\}
\]
Let further $A :\R^n \larrow \R^N$ be an affine map. 

\ms

a) Suppose that for some $\Om'\Subset \Om$ and any $\la >0$, $u$ satisfies
\[
 \|Du\|_{L^\infty(\Om')}\,\leq\, \big\|Du+\la\, DA\big\|_{L^\infty(\Om')}.
\]
Then, we have
\[
 \max_{z\in \overline{\Om'}} \, \big\{Du(z) : DA \big\}\, \geq\, 0.
\]

b) Given $x\in \Om$ and $0<\e <\dist(x,\p\Om)$, the set
\[
 \Om_\e (x)\, :=\, \Big\{ y\in \Om\, \big| \, |Du(y)|<|Du(x)|\Big\} \cap \mB_\e(x)
\]
is open and compactly contained in $\Om$ with $x\in \big(\Om_\e(x)\big)(u)$ if non-empty, that is
\[
 |Du(x)|\, =\, \|Du\|_{L^\infty(  \Om_\e (x) )}.
\]
\end{lemma}

\BPL \ref{lemma1}. a) By assumption we have
\[
\|Du\|^2_{L^\infty(\Om')}\,\leq\, \|Du+\la\,DA\|^2_{L^\infty(\Om')}
\]
and hence
\[
\begin{split}
\underset{\Om'}{\ess\, \sup}\, |Du|^2\, &\leq \, \underset{\Om'}{\ess\, \sup}\,\Big\{|Du|^2 \, +\, 2\la\, Du : DA\, +\, \la^2|DA|^2\Big\}\\
&\leq \, \underset{\Om'}{\ess\, \sup}\, |Du|^2 \, +\, 2\la\, \underset{\Om'}{\ess\, \sup}\,\big\{ Du : DA\big\} \, +\, \la^2|DA|^2.
\end{split}
\]
Consequently,
\[
 \underset{\Om'}{\ess\, \sup}\,\big\{ Du : DA\big\} \, +\, \frac{\la}{2}|DA|^2 \, \geq\, 0
\]
and by letting $\la \ri 0^+$, we obtain the desired inequality. b) is immediate from the definitions.     \qed

\ms

Lemma \ref{lemma1} above is in general true for locally Lipschitz maps, once we replace $|Du|$ by the \emph{local $L^\infty$ norm}
\[
\|Du\|_\infty(x)\, :=\, \lim_{\e\ri0} \|Du\|_{L^\infty(\mB_\e(x))}
\]
which has enough upper semi-continuity properties.

\begin{lemma} \label{lemma2} Let $\Om\sub \R^n$ be open and $u\in C^1(\Om,\R^N)$.  Given $\Om'\Subset \Om$, let $\Om'(u)$ be as in Lemma \ref{lemma1}. Let further $A :\R^n \larrow \R^N$ be an affine map. We set
\[
h(t)\, :=\,  \big\|Du+t\, DA\big\|^2_{L^\infty(\Om')}  - \|Du\|^2_{L^\infty(\Om')},\ \ \ t\geq0.
\]
Then, $h$ is convex, $h(0)=0$ and also the lower right Dini derivative of $h$ at zero satisfies
\[
\underline{D} h(0^+)\, :=\, \underset{t\ri 0^+}{\lim\inf} \, \frac{h(t)-h(0)}{t}  \, \geq \, \max_{y\in {\Om'(u)}} \, \big\{ 2\,Du(y):DA \big\}.
\]
\end{lemma}

\BPL \ref{lemma2}. Effectively, this is an application of Danskin's theorem \cite{D}, but we may also prove it directly. By setting
\[
H(t,y)\, :=\,  \big|Du(y)+t\, DA\big|^2
\]
we have 
\[
h(t)\, =\, \max_{y\in \overline{\Om'}}\, H(t,y) - \max_{y\in \overline{\Om'}} \,H(0,y)
\]
and also for any $t\geq0$ the maximum $\max_{y\in \overline{\Om'}}\, H(t,y)$ is realised at (at least one) point $y^t \in \overline{\Om'}$. Hence 
\[
\begin{split}
\frac{1}{t}\big( h(t)-h(0)\big)\, &=\, \frac{1}{t}\Big[\max_{y\in \overline{\Om'}}\, H(t,y) \,-\, \max_{y\in \overline{\Om'}} \,H(0,y)\Big] \\
&=\, \frac{1}{t}\Big[   H(t,y^t) \, -\, H(0,y^0) \Big] \\
&=\, \frac{1}{t}\Big[ \big( H(t,y^t) -  H(t,y^0) \big)\,+\,\big( H(t,y^0)- H(0,y^0)\big)\Big] \\
&\geq \, \frac{1}{t} \big( H(t,y^0)\,-\, H(0,y^0) \big),
\end{split} 
\]
where $y^0\in \overline{\Om'}$ is any point such that 
\[
|Du(y^0)|\, =\, H(0,y^0)\, =\, \max_{\overline{\Om'}}H(0,\cdot)\, =\, \|Du\|_{L^\infty(\Om')}. 
\]
Hence, by the definition of the set $\Om'(u)$ in Lemma \ref{lemma1}, we have
\[
\begin{split}
\underline{D} h(0^+)  \, &=\ \underset{t\ri 0^+}{\lim\inf} \, \frac{1}{t}\big( h(t)-h(0)\big) 
\\
& \geq \, \max_{y\in {\Om'(u)}} \left\{ \underset{t\ri 0^+}{\lim\inf} \, \frac{1}{t} \Big( H(t,y)- H(0,y) \Big) \right\}\\
& = \, \max_{y\in {\Om'(u)}} \left\{ \underset{t\ri 0^+}{\lim\inf} \, \frac{1}{t} \Big( \big|Du(y)+t\, DA\big|^2- |Du(y)|^2 \Big) \right\}\\
&= \, \max_{y\in {\Om'(u)}} \, \big\{ 2\,Du(y):DA \big\}.
\end{split}
\]
The lemma follows. \qed

\smallskip

Let us also record for later use the elementary inequality
\[
h(t)\,-\,h(0)\, \geq\, \underline{D} h(0^+)\,t,\ \ \ t\geq0,
\]
which is an immediate consequence of the definitions of convexity and of the lower right Dini derivative.

\section{The scalar case $N=1$}

The following is the first main result of this section. 

\begin{theorem}[$C^1$ $\infty$-Harmonic functions] \label{theorem8} Let $\Om\sub \R^n$ be open and $u\in C^1(\Om)$. Given $\Om'\Subset \Om$, let $\Om'(u)$ be as in Lemma \ref{lemma1} and consider the sets of affine functions
\[
\mA^{\pm,\infty}_{\Om'}(u)\, :=\, \left\{A :\, \R^n \ri \R \, \left|
\begin{array}{l}
D^2A\, \equiv\, 0 \text{ and there exist } \xi\in \R^{\pm},\,   \\ 
x\in \Om'(u) \text{ and } \X_x  \in D^{2,\pm}u(x) \\
\text{s.\,th.: } DA\, \equiv \,\xi \, \X_x Du(x)
 \end{array}
 \right. \!\!\!
 \right\} \bigcup \, \R.
\]
Then, if $|D u|$ has no minima in $\Om$, we have the equivalences
\[
\left.
\begin{array}{l}
 Du \ot Du :D^2u \geq  0 \text{ on } \Om, \ms\\
\text{in the Viscosity sense}
\end{array} \right\}
\ \ \Longleftrightarrow\ \ 
\left\{
\begin{array}{l}
\text{For all }\, \Om'\Subset \Om \text{ and }  A\in \mA^{+,\infty}_{\Om'}(u), \ms\\
 \|Du\|_{L^\infty(\Om')}\,\leq\, \|Du+DA\|_{L^\infty(\Om')},
 \end{array}
 \right.
\]
\[
\left.
\begin{array}{l}
 Du \ot Du :D^2u \leq  0 \text{ on }\Om, \ms\\
\text{in the Viscosity sense}
\end{array} \right\}
\ \ \Longleftrightarrow\ \ 
\left\{
\begin{array}{l}
\text{For all }\, \Om'\Subset \Om \text{ and }  A\in \mA^{-,\infty}_{\Om'}(u), \ms\\
 \|Du\|_{L^\infty(\Om')}\,\leq\, \|Du+DA\|_{L^\infty(\Om')}.
 \end{array}
 \right.
\]

\end{theorem}

We note that by the $C^1$ regularity results for $\infty$-Harmonic functions of Savin and Evans-Savin \cite{S,ES}, if $n=2$ the hypothesis that $u$ is a $C^1(\Om)$ viscosity solution is superfluous.

Obviously, for certain subdomains it may happen that $\mA^{\pm,\infty}_{\Om'}(u)$ contain only the trivial (i.e.\ constant) functions if $J^{2,\pm}u(x)=\emptyset$ for all points $x\in \Om'(u)$. Hence, the minimality property above with respect to affine functions is an effective restatement of the definition of viscosity sub/super solutions.

In the event that the solution is smooth, Theorem \ref{theorem8} above simplifies to the following statement for classical solutions of the $\infty$-Laplacian: 

\begin{corollary}[$C^2$ $\infty$-Harmonic functions]  \label{corollary9}  Suppose that $\Om\sub \R^n$ is open and $u\in C^2(\Om)$. Then, if $|D u|$ has no minima in $\Om$ we have the equivalence
\[
Du \ot Du :D^2u \, = \, 0 \text{ on }\Om 
\ \ \Longleftrightarrow\ \ 
\left\{
\begin{array}{l}
\text{For all }\, \Om'\Subset \Om \text{ and }  A\in \big(\mA^{+,\infty}_{\Om'} \cup  \mA^{-,\infty}_{\Om'}\big) (u), \ms\\
 \|Du\|_{L^\infty(\Om')}\,\leq\, \|Du+DA\|_{L^\infty(\Om')}
 \end{array}
 \right.
\]
\[
\hspace{78pt} \Longleftrightarrow\ \ 
\left\{
\begin{array}{l}
\text{For all }\, \Om'\Subset \Om \text{ and }  A\in \mA^{\infty}_{\Om'}(u), \ms\\
 \|Du\|_{L^\infty(\Om')}\,\leq\, \|Du+DA\|_{L^\infty(\Om')}.
 \end{array}
 \right.
\]
Here $\mA^{\infty}_{\Om'}(u)$ is the set of affine functions 
\[
\mA^{\infty}_{\Om'}(u)\, =\, \left\{A :\, \R^n \ri \R \, \left|
\begin{array}{l}
D^2A\, \equiv\, 0 \text{ and there exist } \xi\in \R,\,   \\ 
\text{and }x\in \Om'(u) \text{ s.\,th.\ $A$ is parallel}\\
\text{to the tangent of $\xi|Du|^2$ at }x 
 \end{array}
 \right. \!\!\!
 \right\}.
\]
\end{corollary}

\BPT \ref{theorem8}. Suppose that for any $\Om'\Subset \Om$ and any affine function in $\mA^{+,\infty}_{\Om'}(u)$, we have
\[
\|Du\|_{L^\infty(\Om')}\,\leq\, \|Du+DA\|_{L^\infty(\Om')}.
\]
Fix any $x\in \Om$ such that $(Du(x),\X_x)\in J^{2,+}u(x)$, whence $\X_x\in D^{2,+}u(x)$. Consider the affine function 
\[
A(z)\, :=\, \xi\, \X_x : Du(x) \ot (z-x),\ \ \ z\in \R^n,
\]
where $\xi \geq0$. Fix also $\e>0$ and let $\Om_\e(x)$ be as in Lemma \ref{lemma1}. Then,  for any $\la>0$, the affine function $\la A$ is contained in $\mA^{+,\infty}_{\Om_\e(x)}(u)$. Hence,
\[
\|Du\|_{L^\infty(\Om_\e(x))}\,\leq\, \|Du+\la\,DA\|_{L^\infty(\Om_\e(x))}.
\]
By applying Lemma \ref{lemma1} to $u$ and $A$, we have
\[
\begin{split}
0 \,  &\leq\, \max_{z\in \overline{\Om_\e(x)}} \, \big\{ Du(z) \cdot DA \big\} \\
&=\, \max_{z\in \overline{\Om_\e(x)}} \big\{ Du(z) \cdot\big(\xi\, \X_x :Du(x)) \big\}  \\
&=\, \max_{z\in \overline{\Om_\e(x)}} \big\{ \xi\, \big(\X_x :Du(x) \ot Du(z) \big)\big\} \\
&\!\! \larrow\, \xi\, \big(\X_x :Du(x) \ot Du(x) \big), 
\end{split}
\] 
as $\e \ri 0$. Hence, $Du \ot Du :D^2u \geq 0$ on $\Om$ in the viscosity sense. 

\smallskip

Conversely, fix any $\Om'\Subset \Om$ and $x\in \Om'(u)$. If it happens $J^{2,+}u(x)\neq \emptyset$, then any $A\in \mA^{+,\infty}_{\Om'}(u)$ can be written as 
\[
A(z)\, =\, a\, +\, \xi\, \X_x : Du(x) \ot z ,\ \ \ z\in \R^n,
\]
for some $a\in\R$, $\xi \geq0$ and $\X_x \in D^{2,+}u(x)$. Let $h$ be the function of Lemma \ref{lemma2} for such an $A$.  By applying Lemma \ref{lemma2} to this setting, we have
\[
\begin{split}
\underline{D}h(0^+)\, &\geq\, \max_{y\in {\Om'(u)}} \, \big\{ 2\,Du(y) \cdot DA \big\} \\
& \geq \, 2\,Du(x) \cdot DA \\
&=\, 2\, Du(x) \cdot\big(\xi\, \X_x :Du(x)) \big\}  \\
& =\, 2\,\xi\, \big(\X_x :Du(x) \ot Du(x) \big)\\
& \geq \, 0, 
\end{split}
\]
since by assumption $Du \ot Du :D^2u \geq 0$ on $\Om$ in the viscosity sense. Since $h(0)=0$ and $h$ is convex, it follows that 
\[
h(t) \, \geq\, h(0)\, +\, \underline{D}h(0^+)\,t \, \geq \, 0, \ \ \  t\geq0,
\]
and hence, by the definition of $h$ we obtain
\[
\|Du\|_{L^\infty(\Om')}\,\leq\, \|Du+DA\|_{L^\infty(\Om')}
\]
for any $\Om'\Subset \Om$ and any $A\in \mA^{+,\infty}_{\Om'}(u)$. The case of supersolutions follows similarly and hence the theorem has been established.   \qed

\ms

\BPCOR \ref{corollary9}. The first equivalence of the statement is immediate. Since by assumption $u\in C^2(\Om)$, we have that 
\[
 J^{2,+}u(x) \cap J^{2,-}u(x)\, =\, \big\{ \big(Du(x),D^2u(x)\big) \big\}
\]
and hence $D^{2,+}u(x) \cap D^{2,-}u(x) = \{D^2u(x)\}$. The second equivalence of the statement follows by making the choice $\X_x \in D^{2,\pm}u(x)$ in the proof of Theorem \ref{theorem8} above and repeating all the steps. Then, by noting that
\[
\X_x Du(x)\, =\, D\left( \frac{1}{2}|Du|^2\right)(x)
\] 
it follows that for any $\Om'\Subset \Om$ the set $\mA^{\infty}_{\Om'}(u)$ contains only affine functions of the form 
\[
A(z)\, =\, a\,+\, \xi\, D\big(|Du|^2\big)(x) \cdot (z-x), \ \ \ z\in \R^n,
\]
for $a,\xi \in \R$ and $x\in \Om'(u)$. The corollary ensues.        \qed

\ms

Theorem \ref{theorem8} extends relatively easily to the case of the $p$-Laplacian for $1<p<\infty$ which, quite surprisingly, can also be characterised by the $L^\infty$ functional via affine variations. In view of the well known $C^{1,\al}$ regularity results for $p$-Harmonic mappings \cite{U}, the hypothesis that solutions are $C^1$ is actually superfluous. 

\begin{theorem}[$p$-Harmonic functions]  \label{theorem10} Let $\Om\sub \R^n$ be open and $u\in C^1(\Om)$. Given $\Om'\Subset \Om$, let $\Om'(u)$ be as above and consider the sets of affine functions
\[
\mA^{\pm,p}_{\Om'}(u)\, :=\, \left\{A :\, \R^n \ri \R \, 
\left| 
\begin{array}{l}
D^2A\, \equiv\, 0 \text{ and there exist } \xi\in \R^{\pm},  \\ 
 x\in \Om'(u) \text{ and } \X_x  \in D^{2,\pm}u(x) \text{ s.\,th.: } \\
 DA  \equiv \xi  \Big( \!(p-2)\X_x +(I\!:\!\X_x) I\Big)Du(x)\\
 \end{array}
 \right.\!\!\!\!
 \right\} \bigcup \, \R,
\]
where $p\in (1,\infty)$. Then, if $|D u|$ has no minima in $\Om$ the following statements are equivalent:

\noi (a) We have
\[
\text{div}\big( |Du|^{p-2}Du\big)\, \geq \, 0, \ \text{ weakly on }\Om.
\]
\noi (b) We have
\[
 \left((p-2)Du \ot Du  +  |Du|^2 I\right) : D^2u \, \geq \, 0 
\text{ on } \Om,
\]
in the feeble Viscosity sense.

\noi (c) For all $\Om'\Subset \Om$ and all $A\in \mA^{+,p}_{\Om'}(u) $, we have
\[
 \|Du\|_{L^\infty(\Om')}\,\leq\, \|Du+DA\|_{L^\infty(\Om')}.
\]
The case ``$\,\leq 0$" of supersolutions is symmetrical and corresponds to $\mA^{-,p}_{\Om'}(u)$ as in Theorem \ref{theorem8} above.
\end{theorem}

In the case of the usual Laplacian for $p=2$, the affine functions in $\mA^{+,2}_{\Om'}(u)$ of Theorem \ref{theorem10} satisfy $DA=\xi (\X_x\!:\! I)Du(x)$, where $\xi\geq 0$, $\X_x  \in D^{2,\pm}u(x)$,  $\Om'\Subset \Om$ and $x\in \Om'(u)$.

\BPT \ref{theorem10}. The idea is similar to that of the proof of Theorem  \ref{theorem8}, so we basically need to indicate the points where it differs. We begin by noting by the results of the papers \cite{K0,JLM,JJ}, it follows that a function is weakly $p$-subharmonic on $\Om$ (that is we have $\text{div}\big( |Du|^{p-2}Du\big)\geq 0$ holding weakly on $\Om$) if and only if it is $p$-subharmonic on $\Om$ in the feeble viscosity sense for the $p$-Laplacian expanded: 
\[
|Du|^{p-4}\Big((p-2)Du \ot Du  +  |Du|^2 I\Big):D^2u \geq \, 0, \ \ \text{ on }\Om.
\]
Since by definition of the feeble Jets we do not check anything in the viscosity criterion when the gradient vanishes, the $p$-Laplacian is equivalent in the feeble viscosity sense to
\[
\Big((p-2)Du \ot Du  +  |Du|^2 I\Big):D^2u \geq \, 0, \ \ \text{ on }\Om.
\]
As a consequence, $(a) \Leftrightarrow (b)$. We suppose now that for any $\Om'\Subset \Om$ and any affine function $A\in \mA^{+,\infty}_{\Om'}(u)$, we have
\[
\|Du\|_{L^\infty(\Om')}\,\leq\, \|Du+DA\|_{L^\infty(\Om')}.
\]
Fix any $x\in \Om$ such that $(Du(x),\X_x)\in J_0^{2,+}u(x)$, whence $\X_x\in D^{2,+}u(x)$. Consider the affine function 
\[
A(z)\, :=\, \Big( (p-2)\X_x\, +\, (I: \X_x ) I\Big): Du(x) \ot (z-x),\ \ \ z\in \R^n.
\]
Fix also $\e>0$ and let $\Om_\e(x)$ be as in Lemma \ref{lemma1} and note that  for any $\la>0$, $\la A \in \mA^{+,p}_{\Om_\e(x)}(u)$. Hence, by arguing as in Theorem  \ref{theorem8} we have that
\[
\begin{split}
0\, &\leq\,  Du(x) \cdot DA\\
&=\, Du(x) \cdot  \Big( (p-2)\X_xDu(x)\, +\, (I: \X_x ) Du(x)\Big) \\
&=\,\Big( (p-2) Du(x) \ot Du(x)\, +\, |Du(x)|^2I \Big): \X_x.
\end{split}
\] 
Hence, $u$ is a feeble viscosity solution on $\Om$. 

\smallskip

Conversely, fix any $\Om'\Subset \Om$ and $x\in \Om'(u)$. If $J_0^{2,+}u(x)\neq \emptyset$, then any $A\in \mA^{+,p}_{\Om'}(u)$ can be written as 
\[
A(z)\, =\, a\, +\, \xi\, \Big( (p-2)\X_x\, +\, (I: \X_x ) I\Big): Du(x) \ot z ,\ \ \ z\in \R^n,
\]
for some $a\in\R$, $\xi \geq 0$ and some $(Du(x),\X_x)\in J_0^{2,+}u(x)$. Let $h$ be the function of Lemma \ref{lemma2} for such an $A$.  By applying Lemma \ref{lemma2}, we have
\[
\begin{split}
\underline{D}h(0^+)\,  & \geq \, 2\,Du(x) \cdot DA \\
& =\, 2\,\xi\, \Big(  (p-2)Du(x) \ot Du(x) : \X_x\, +\, |Du(x)|^2 I: \X_x  \Big)\\
& \geq \, 0, 
\end{split}
\]
since by assumption $u$ is a subsolution on $\Om$ in the feeble viscosity sense. By using that $h(0)=0$ and that $h$ is convex, we deduce as in Theorem  \ref{theorem8} that $h(t)\geq0$ for  $ t\geq0$
and hence
\[
\|Du\|_{L^\infty(\Om')}\,\leq\, \|Du+DA\|_{L^\infty(\Om')}
\]
for any $A\in \mA^{+,p}_{\Om'}(u)$ and any $\Om'\Subset \Om$. Thus, $(b) \Leftrightarrow (c)$. The case of supersolutions follows analogously and hence the theorem ensues.  \qed

\section{The vectorial case $N\geq 2$}

In this section we extend the results of the previous section to the full case of the $\infty$-Laplace system. We begin by noting that \eqref{1.1} actually consists of two independent systems, the second of which is identically trivial in the scalar case. Namely, if $u:\Om\sub \R^n \larrow \R^N$ is smooth, then 
\[
\De_\infty u =0 \ \ \Longleftrightarrow\ \ \left\{
\begin{split}
Du \ot Du :D^2u\, &=\, 0,\\
|Du|^2[Du]^\bot \De u\, &=\,0.
\end{split}
\right.
\]
This is an immediate consequence of the mutual perpendicularity of the vector fields $Du \ot Du :D^2u$ and $|Du|^2[Du]^\bot \De u$; indeed, it suffices to recall that $[Du]^\bot$ is the projection on the orthogonal complement of $R(Du)$ and to note the identity
\[
2\, Du \ot Du :D^2u \, =\,  Du\, D\big(|Du|^2 \big).
\]
Our last main result is the following:

\begin{theorem}[$C^1$ $\infty$-Harmonic mappings]  \label{theorem11} Let $\Om\sub \R^n$ be open and $u\in C^1(\Om,\R^N)$ with $|D u|$ having no minima in $\Om$. Given a set $\Om'\Subset \Om$, let $\Om'(u)$ be as in Lemma \ref{lemma1}. Consider first the set of affine maps
\[
\mA^{\top,\infty}_{\Om'}(u) := \left\{ \! A :\, \R^n \ri \R^N \,  
\left|
\begin{array}{l}
D^2A\, \equiv\, 0\ \& \text{ there exist } \xi\in \R^N, \, x\in \Om'(u) \\ 
\mD^2u \in \mY\big(\Om,   \smash{\overline{\R}}^{Nn^2}_s \big)\ \&\ \X_x  \in \supp_*\big(\mD^2u(x)\big) \\
\text{s.\,th.: }  DA\, \equiv \,\xi \ot \big(\X_x : Du(x)\big)
 \end{array}
 \right. \!\!\!\!\!
 \right\} \bigcup \, \R^N.
\]
Then, we have the equivalence
\[
\left.
\begin{array}{l}
Du \ot Du :D^2u \,=\, 0  \ms\\
\text{on } \Om, \text{ in the $\mD$-sense}
\end{array} \right\}
\ \ \Longleftrightarrow\ \ 
\left\{
\begin{array}{l}
\text{For all }\, \Om'\Subset \Om \text{ and }  A\in \mA^{\top,\infty}_{\Om'}(u) , \ms\\
 \|Du\|_{L^\infty(\Om')}\,\leq\, \|Du+DA\|_{L^\infty(\Om')}.
 \end{array}
 \right.
\]
Further, consider the set of affine maps
\[
\mA^{\bot,\infty}_{\Om'}(u) := \left\{ \! A :\, \R^n \ri \R^N \, 
\left|  
\begin{array}{l}
D^2A\, \equiv\, 0\ \& \text{ there exist } x\in \Om'(u),\, \mD^2u \in \\ 
 \mY\big(\Om,   \smash{\overline{\R}}^{Nn^2}_s \big)\ \&\, \X_x  \in \supp_*\big(\mD^2u(x)\big) \text{ s.\,th.:}\\
 A(x) \in R\big(Du(x)\big)^\bot\ \&\ DA \in \mathscr{L}^{\X_x}\big(A(x)\big)
 \end{array}
 \right. \!\!\!\!\!
 \right\} \bigcup \, \R^N
\]
where for any $a\in \R^N$, $\mathscr{L}^{\X_x}(a)$ is an affine matrix space defined as
\[
\mathscr{L}^{\X_x}(a)\, :=\, 
\left\{
\begin{array}{l}
\Big\{ X\in \R^{Nn}\, \Big| \, Du(x):X = -\big(a\ot I\big)\!:\!\X_x\Big\}, \text{ if }Du(x)\neq0 \ms\\
\{0\}, \hspace{163pt} \text{ if }Du(x)=0.
\end{array}
\right.
\]
Then, we have the equivalence
\[
\left.
\begin{array}{l}
|Du|^2[Du]^\bot \De u \,=\, 0 \ms\\
\text{on } \Om, \text{ in the $\mD$-sense}
\end{array} \right\}
\ \ \Longleftrightarrow\ \ 
\left\{
\begin{array}{l}
\text{For all }\, \Om'\Subset \Om \text{ and }  A\in \mA^{\bot,\infty}_{\Om'}(u), \ms\\
 \|Du\|_{L^\infty(\Om')}\,\leq\, \|Du+DA\|_{L^\infty(\Om')}.
 \end{array}
 \right.
\]
\end{theorem}

In view of Theorem \ref{theorem11} above, a mapping is $\infty$-Harmonic in the $\mD$-sense iff it minimises with respect to the union of the sets of affine variations of the tangential and the normal component:
\[
\left.
\begin{array}{l}
\De_\infty u \,=\, 0 \text{ on } \Om,  \ms\\
 \text{in the $\mD$-sense}
\end{array} \right\}
\ \ \Longleftrightarrow\ \ 
\left\{
\begin{array}{l}
\text{For all }\, \Om'\Subset \Om \text{ and }  A\in \big(\mA^{\top,\infty}_{\Om'} \cup  \mA^{\bot,\infty}_{\Om'} \big) (u) , \ms\\
 \|Du\|_{L^\infty(\Om')}\,\leq\, \|Du+DA\|_{L^\infty(\Om')}.
 \end{array}
 \right.
\]
In the event that $u\in C^2(\Om,\R^N)$, Theorem \ref{theorem11} simplifies to the following statement for classical solutions of the $\infty$-Laplace system:

\begin{corollary}[$C^2$ $\infty$-Harmonic mappings]  \label{corollary12}  Suppose that $\Om\sub \R^n$ is open and $u\in C^2(\Om,\R^N)$. Then, if $|D u|$ has no minima in $\Om$ we have the equivalence
\[
\De_\infty u \, = \, 0 \text{ on }\Om 
\ \ \Longleftrightarrow\ \ 
\left\{
\begin{array}{l}
\text{For all }\, \Om'\Subset \Om \text{ and }  A\in \big(\mA^{\top,\infty}_{\Om'}\cup \mA^{\bot,\infty}_{\Om'}\big)(u), \ms\\
 \|Du\|_{L^\infty(\Om')}\,\leq\, \|Du+DA\|_{L^\infty(\Om')},
 \end{array}
 \right.
\]
where $\mA^{\top,\infty}_{\Om'}(u)$, $\mA^{\bot,\infty}_{\Om'}(u)$ are the sets of affine maps
\[
\mA^{\top,\infty}_{\Om'}(u)\, =\, \left\{A :\, \R^n \ri \R^N \, \left|
\begin{array}{l}
D^2A\, \equiv\, 0 \text{ and there exist } \xi\in \R^N,\,   \\ 
\text{and }x\in \Om'(u) \text{ s.\,th.\ $A$ is parallel}\\
\text{to the tangent of $\xi|Du|^2$ at }x 
 \end{array}
 \right. \!\!\!
 \right\}, \ \ \ 
\]
\[
\mA^{\bot,\infty}_{\Om'}(u)\, =\, \left\{A :\, \R^n \ri \R^N \, \left|
\begin{array}{l}
D^2A\, \equiv\, 0 \text{ and there exists }x\in \Om'(u)   \\ 
\text{such that $A$ is normal to $Du$ at $x$ }\\
\text{and $A^\top \!Du$ is divergenceless at $x$ } 
 \end{array}
 \right. \!\!\!
 \right\}.
\]
\end{corollary}

\BPT \ref{theorem11}. We begin by a general observation about the notion of $\mD$-solutions $u:\Om\sub \R^n \larrow \R^N$ in $C^1(\Om,\R^N)$ to a homogeneous 2nd order quasilinear system of the form
\[
\A(Du): D^2u\, =\, 0, \ \ \text{ on }\Om,
\] 
when $\A$ is Borel measurable. By definition, every diffuse hessian $\mD^2u \in \mY\big(\Om,   \smash{\overline{\R}}^{Nn^2}_s \big)$ of a candidate solution $u$ is defined a.e.\ on $\Om$ as a weakly* measurable probability valued map $\Om \larrow  \smash{{\R}}^{Nn^2}_s\cup\{\infty\}$. Hence, we may modify each $\mD^2u$ on a Lebesgue nullset and choose from each equivalence class the representative which is redefined as $\de_{\{0\}}$ at points where $\mD^2u(x)$ does not exist. Moreover, let $u$ be a fix map in $C^1(\Om,\R^N)$. Since $Du(x)$ exists for all $x\in \Om$, by perhaps a further re-definition of every $\mD^2u$ on a Lebesgue nullset, it follows that $u$ is $\mD$-solution to the system if and only if for (any fixed such representative of) any diffuse hessian, we have
\[
\A\big(Du(x)\big): \X_x\, =\, 0, \ \ \text{ for all }x\in\Om \text{ and } \X_x \in \supp_*\big(\mD^2u(x)\big).
\] 
(We remind that at points $x\in \Om$ for which  $\mD^2u(x) =\de_{\{\infty\}}$ and hence $\supp_*\big(\mD^2u(x)\big)$ $= \emptyset$, the above condition is understood as being trivially satisfied.) We will apply this observation to the two independent systems 
\[
Du\ot Du :D^2u\,=\,0,\ \ \ \ |Du|^2\big([Du]^\bot \!\ot I \big):D^2u\,=\,0
\]
comprising the $\infty$-Laplace system.

Suppose now that for some $\Om'\Subset \Om$ and some affine mapping $A\in \mA^{\top,\infty}_{\Om'}(u)$, we have
\[
\|Du\|_{L^\infty(\Om')}\,\leq\, \|Du+DA\|_{L^\infty(\Om')}.
\]
Fix any $x\in \Om$ and any diffuse hessian $\mD^2u \in \mY\big(\Om,   \smash{\overline{\R}}^{Nn^2}_s \big)$ such that $\supp_*\big(\mD^2u(x)\big)$ 
$\neq \emptyset$ and pick any $\X_x \in \supp_*\big(\mD^2u(x)\big)$. Fix also $\xi \in \R^N$ and consider the affine map which is defined by
\[
A(z)\, :=\, \xi \ot \big(\X_x : Du(x)\big) \cdot (z-x),\ \ \ z\in \R^n.
\]
In index form this means
\[
A_\al(z)\, =\, \xi_\al \sum_{\be=1}^N\sum_{i,j=1}^n\Big((\X_x)_{\be ji}  D_ju_\be (x)\Big) (z-x)_i,\ \ \ \al=1,..,N.
\]
For $\e>0$ small, let $\Om_\e(x)$ be as in Lemma \ref{lemma1}. Then,  $\la A \in \mA^{\top,\infty}_{\Om_\e(x)}(u)$ for any $\la >0$. Thus,
\[
\|Du\|_{L^\infty(\Om_\e(x))}\,\leq\, \|Du+\la\,DA\|_{L^\infty(\Om_\e(x))}
\]
and by applying Lemma \ref{lemma1} to $u$ and $A$, we have
\[
\begin{split}
0 \, &\leq\, \max_{z\in \overline{\Om_\e(x)}} \Big\{ Du(z) : \Big(\xi \ot \X_x :Du(x)\Big) \Big\} \\
& =\, \max_{z\in \overline{\Om_\e(x)}} \left\{ \sum_{\al=1}^N\sum_{i=1}^n \, D_i u_\al(z)\, \xi_\al \sum_{\be=1}^N\sum_{j=1}^n \, (\X_x)_{\be ji}  D_ju_\be (x)  \right\} 
\\
&\leq\,  \max_{z\in \overline{\Om_\e(x)}} \left\{ \sum_{\al,\be=1}^N\sum_{i,j=1}^n \, \xi_\al \, D_i u_\al(z)\, D_ju_\be (x) \,  (\X_x)_{\be ji}  \right\} \\
 &\!\! \larrow\, \sum_{\al,\be=1}^N\sum_{i,j=1}^n \, \xi_\al \, \Big(D_i u_\al(x) \,  D_ju_\be (x) \, (\X_x)_{\be ji}\Big)
\end{split}
\] 
as $\e \ri 0$, and hence
\[
 \xi \cdot \big(Du(x) \ot Du(x) :\X_x\big) \, \geq \, 0,
\] 
for any $\xi \in \R^N$. By the arbitrariness of $\xi$ we deduce that $Du(x) \ot Du(x) :\X_x=0$. As a consequence, $Du \ot Du :D^2u = 0$ in the $\mD$-sense on $\Om$.

Now we argue similarly for the normal component of the system. Suppose that for any $\Om'\Subset \Om$ and any $A\in \mA^{\bot,\infty}_{\Om'}(u)$, we have
\[
\|Du\|_{L^\infty(\Om')}\,\leq\, \|Du+DA\|_{L^\infty(\Om')}.
\]
We fix as before $x\in \Om$ and $\X_x \in \supp_*\big(\mD^2u(x)\big)$. If $Du(x)=0$, then the system $|Du|^2[Du]^\bot \De u=0$ is trivially satisfied at $x$. If $Du(x)\neq 0$, then we choose any direction normal to $Du(x)$, that is 
\[
n_x\, \in R\big(Du(x)\big)^\bot \sub \,\R^N, 
\]
which means that $n_x^\top Du(x)=0$. We note that if $Du(x) : \R^n \larrow \R^N$ is surjective, then we can find only the trivial $n_x=0$, but the system  $|Du|^2[Du]^\bot \De u=0$ is satisfied at $x$ anyhow because $[Du(x)]^\bot =0$. We also fix any matrix $N_x$ in the affine space $\mathscr{L}^{\X_x}(n_x)$. By the definition of $\mathscr{L}^{\X_x}(n_x)$, this means that
\[
N_x : Du(x)\, =\, -(n_x \ot I):\X_x.
\]
We consider the affine map which is defined by
\[
A(z)\, :=\, n_x\, +\, N_x (z-x),\ \ \ z\in \R^n.
\]
We now claim that $\la A\in \mA^{\bot,\infty}_{\Om'}(u)$ for any $\la \in\R$. Indeed, this is a consequence of our choices and of the following homogeneity property of the space $\mathscr{L}^{\X_x}(a)$: 
\[
\mathscr{L}^{\X_x}(\la a) \, =\, \la\, \mathscr{L}^{\X_x}(a), \ \ \ \la \in \R. 
\] 
Hence, we have
\[
\|Du\|_{L^\infty(\Om')}\,\leq\, \|Du+\la\, DA\|_{L^\infty(\Om')}.
\]
By applying Lemma \ref{lemma1} to $u$ and $A$, we have
\[
\begin{split}
0 \, &\leq\, \max_{z\in \overline{\Om_\e(x)}} \big\{ Du(z) : N_x \big\}  \\
& \!\!\larrow\, Du(x) : N_x \\
& =\, -(n_x \ot I):\X_x, 
\end{split}
\] 
as $\e \ri 0$. Hence, we have $(n_x \ot I):\X_x\leq 0$ and by the arbitrariness of the direction $n_x \, \bot\, R\big(Du(x)\big)$, we obtain that $(n_x \ot I):\X_x=0$. Thus, $\big([Du(x)]^\bot \ot I \big):\X_x =0$ and as a consequence $|Du|^2[Du]^\bot \De u=0$ in the $\mD$-sense on $\Om$.

\smallskip

Conversely, we fix $\Om'\Subset \Om$ and $x\in \Om'(u)$ and any $A\in \mA^{\top,\infty}_{\Om'}(u)$ corresponding to a diffuse hessian $\mD^2u \in \mY\big(\Om, \smash{\overline{\R}}^{Nn^2}_s \big)$ and some $\X_x \in \supp_*(\mD^2u(x))$ and $\xi \in \R^N$. We take as $h$ to be the function of Lemma \ref{lemma2}. By applying Lemma \ref{lemma2} to this setting, we have
\[
\begin{split}
\underline{D}h(0^+)\, &\geq\, \max_{y\in {\Om'(u)}} \, \big\{ 2\,Du(y) : DA \big\} \\
& \geq \, 2\,Du(x) : DA 
\\
&\geq\, 2 \sum_{\al,\be=1}^N\sum_{i,j=1}^n \, D_i u_\al(x)\, \xi_\al \, (\X_x)_{\be ji}  D_ju_\be (x) 
\end{split}
\]
and hence
\[
\begin{split}
\underline{D}h(0^+)\, &\geq\, 2\xi \cdot \big(Du(x) \ot Du(x) :\X_x\big) \\
& = \, 0, 
\end{split}
\]
since by assumption $Du \ot Du :D^2u= 0$ on $\Om$ in the $\mD$-sense. In view of the fact that $h(0)=0$ and $h$ is convex, it follows that 
\[
h(t) \, \geq\, h(0)\, +\, \underline{D}h(0^+)\,t \, \geq \, 0, \ \ \  t\geq0,
\]
and hence
\[
\|Du\|_{L^\infty(\Om')}\,\leq\, \|Du+DA\|_{L^\infty(\Om')}, \ \ \ A\in \mA^{\top,\infty}_{\Om'}(u),\ \Om'\Subset \Om.
\]
The case of $A\in \mA^{\bot,\infty}_{\Om'}$ is completely analogous: any such nonconstant $A$ satisfies $A(x) \, \bot\, R(Du(x))$ and $DA \in \mathscr{L}^{\X_x}\big(A(x)\big)$ for some $\X_x \in \supp_*(\mD^2u(x))$ and some $x\in \Om'(u)$. By applying Lemma \ref{lemma2} again, we have
\[
\begin{split}
\underline{D}h(0^+)\, &\geq\, \max_{y\in {\Om'(u)}} \, \big\{ 2\,Du(y) : DA \big\} \\
& \geq \, 2\,Du(x) : DA .
\end{split}
\]
If $Du(x)\neq 0$, then by the definition of $\mathscr{L}^{\X_x}\big(A(x)\big)$ we have
\[
\begin{split}
\underline{D}h(0^+)\, &\geq \, 2\, DA : Du(x)\\ 
&=\, -2\, (n_x \ot I):\X_x \\
&=\, - 2\, n_x ^\top \Big(\big( [Du(x)]^\bot \ot I \big) : \X_x \Big) \\
&=\, 0
\end{split}
\]
because by assumption $|Du|^2[Du]^\bot \De u= 0$ on $\Om$ in the $\mD$-sense. If $Du(x)=0$, then again $\underline{D}h(0^+) \geq 0$. In either cases, we obtain
 \[
h(t) \, \geq\, h(0)\, +\, \underline{D}h(0^+)\,t \, \geq \, 0, \ \ \  t\geq0,
\]
and hence
\[
\|Du\|_{L^\infty(\Om')}\,\leq\, \|Du+DA\|_{L^\infty(\Om')}, \ \ \ A\in \mA^{\bot,\infty}_{\Om'}(u),\ \Om'\Subset \Om.
\]
The theorem has been established.   \qed

\ms

\BPCOR \ref{corollary12}. If $u\in C^2(\Om,\R^N)$, then it is an immediate consequence of Lemma \ref{lemma0} that any diffuse hessian of $u$ satisfies
\[
\mD^2u(x)\, =\, \de_{D^2u(x)},\ \ \ x\in\Om,
\]
and by the remarks in the beginning of the proof of Theorem \ref{theorem11}, this happens for all $x\in \Om$. Hence, the only possible $\X_x$ in the reduced support of $\mD^2u(x)$ is $\X_x = D^2u(x)$. For $\mA^{\top,\infty}_{\Om'}$, we have that any possible $A$ satisfies $DA\equiv D\big( \xi |Du|^2)(x)$. For $\mA^{\bot,\infty}_{\Om'}$, we have that any possible $A$ satisfies 
\[
A(x)^\top Du(x)\, =\, 0, \ \ \ \ DA\in \mathscr{L}^{D^2u(x)}\big(A(x)\big), 
\]
which gives
\[
DA : Du(x)\, =\, -\big(A(x) \ot I\big): D^2u(x)\, =\, -A(x) \cdot \De u(x).
\]
Thus, 
\[
\text{div} \big( A^\top Du \big)(x) \, =\,  DA : Du(x)\,  +\, A(x) \cdot \De u(x)\, =\, 0.
\]
The corollary has been established.        \qed

\ms
\ms

\noi {\bf Acknowledgement.} The author would like to thank Craig Evans, Robert Jensen, Jan Kristensen and Juan Manfredi for inspiring scientific discussion relevant to the content of this paper, as well as for their encouragement. He is also indebted to the anonymous referee for the careful reading of the manuscript and for preparing their report so swiftly.

\ms

\bibliographystyle{amsplain}

\begin{thebibliography}{30}

\bibitem[AK]{AK} H. Abugirda, N. Katzourakis, \emph{Existence of $1D$ Vectorial Absolute Minimisers in $L^\infty$ under Minimal Assumptions}, Proceedings of the AMS, accepted. 


\bibitem[AyK]{AyK} B. Ayanbayev, N. Katzourakis, \emph{A Pointwise Characterisation of the PDE System of Vectorial Calculus of Variations in $L^\infty$}, \url{https://arxiv.org/pdf/1611.05936.pdf}.


\bibitem[A1]{A1} G. Aronsson, \emph{Extension of functions satisfying Lipschitz conditions}, Arkiv f\"ur Mat. 6 (1967), 551 - 561.

\bibitem[A2]{A2} G. Aronsson, \emph{On the partial differential equation $u_x^2 u_{xx} + 2u_x u_y u_{xy} + u_y^2 u_{yy} = 0$}, Arkiv f\"ur Mat. 7 (1968), 395 - 425.

\bibitem[BJW1]{BJW1} E. N. Barron, R. Jensen and C. Wang, \emph{The Euler equation and absolute minimizers of $L^{\infty}$ functionals}, Arch. Rational Mech. Analysis 157 (2001), 255 - 283.

\bibitem[BJW2]{BJW2} E. N. Barron, R. Jensen, C. Wang, \emph{Lower Semicontinuity of $L^{\infty}$ Functionals} Ann. I. H. Poincar\'e AN 18, 4 (2001)
495 - 517.

\bibitem[CFV]{CFV} C. Castaing, P. R. de Fitte, M. Valadier, \emph{Young Measures on Topological spaces with Applications in Control Theory and Probability Theory}, Mathematics and Its Applications, Kluwer Academic Publishers, 2004.

\bibitem[C]{C} M. G. Crandall, \emph{A visit with the $\infty$-Laplacian}, in \emph{Calculus of Variations and Non-Linear Partial Differential Equations}, Springer Lecture notes in Mathematics 1927, CIME, Cetraro Italy 2005.

\bibitem[CIL]{CIL} M. G. Crandall, H. Ishii, P.-L. Lions, \emph{User's Guide to
Viscosity Solutions of 2nd Order Partial Differential Equations},
Bulletin of the AMS 27, 1-67 (1992).

\bibitem[CKP]{CKP} G. Croce, N. Katzourakis, G. Pisante, \emph{$\mD$-solutions to the system of vectorial Calculus of Variations in $L^\infty$ via the Baire Category method for the singular values}, ArXiv preprint, \url{http://arxiv.org/pdf/1604.04385.pdf}. 

\bibitem[D]{D} J.M. Danskin, \emph{The theory of min-max with application}, SIAM Journal on Applied Mathematics, 14 (1966), 641 - 664.

\bibitem[Ed]{Ed} R.E. Edwards, \emph{Functional Analysis: Theory and Applications}, Dover Books on Mathematics,  2003.

\bibitem[E]{E} L.C. Evans, \emph{Weak convergence methods for nonlinear partial differential equations}, Regional conference series in mathematics 74, AMS,  1990.

\bibitem[E2]{E2} L.C. Evans, \emph{Partial Differential Equations}, AMS, Graduate Studies in Mathematics Vol. 19, 1998.

\bibitem[ES]{ES} L. C. Evans, O. Savin, \emph{$C^{1,\alpha}$ Regularity for
Infinity Harmonic Functions in Two Dimensions}, Calc. Var. 32, 325 -
347, (2008).

\bibitem[FG]{FG} L.C. Florescu, C. Godet-Thobie, \emph{Young measures and compactness in metric spaces}, De Gruyter, 2012.

\bibitem[FL]{FL} I. Fonseca, G. Leoni, \emph{Modern methods in the Calculus of Variations: $L^p$ spaces}, Springer Monographs in Mathematics, 2007.

\bibitem[JJ]{JJ} V. Julin, P. Juutinen, \emph{A new proof for the equivalence of weak and viscosity solutions for the p-Laplace equation}, Communications in PDE, 37, No 5, 934 - 946, 2012.

\bibitem[JLM]{JLM} P. Juutinen, P. Lindqvist, J.J. Manfedi, \emph{On the equivalence of viscosity solutions and weak solutions for a quasilinear equation}, SIAM J. Math. Anal., Vol. 33, 699 - 717, 2001.

\bibitem[K1]{K1} N. Katzourakis,  \emph{$L^\infty$-Variational Problems for Maps and the Aronsson PDE system}, J.\ Differential Equations, Volume 253, Issue 7 (2012), 2123 - 2139.

\bibitem[K2]{K2} N. Katzourakis,  \emph{Explicit $2D$ $\infty$-Harmonic Maps whose Interfaces have Junctions and Corners}, Comptes Rendus Acad. Sci. Paris, Ser.I, 351 (2013) 677 - 680.

\bibitem[K3]{K3} N. Katzourakis,  \emph{On the Structure of $\infty$-Harmonic Maps}, Communications in PDE, Volume 39, Issue 11 (2014), 2091 - 2124.

\bibitem[K4]{K4} N. Katzourakis, \emph{$\infty$-Minimal Submanifolds}, Proceedings of the AMS, 142 (2014) 2797-2811.

\bibitem[K5]{K0} \emph{Nonsmooth Convex Functionals and Feeble Viscosity Solutions of Singular Euler-Lagrange Equations}, Calculus of Variations and PDE, published online November 2014, (DOI) 10.1007/s00526-014-0786-x.

\bibitem[K6]{K5} N. Katzourakis,  \emph{Nonuniqueness in Vector-valued Calculus of Variations in $L^\infty$ and some Linear Elliptic Systems}, Comm. on Pure and Applied Anal.,  Vol. 14, 1, 313 - 327 (2015). 

\bibitem[K7]{K6} N. Katzourakis,   \emph{Optimal $\infty$-Quasiconformal Immersions},  ESAIM Control Optim. Calc. Var. 21 (2015), no. 2, 561 - 582. 

\bibitem[K8]{K7} N. Katzourakis, \emph{An Introduction to Viscosity Solutions for Fully Nonlinear PDE with Applications to Calculus of Variations in $L^\infty$}, Springer Briefs in Mathematics, 2015, DOI 10.1007/978-3-319-12829-0.

\bibitem[K9]{K8} N. Katzourakis,  \emph{Generalised solutions for fully nonlinear PDE systems and existence-uniqueness theorems}, ArXiv preprint, \url{http://arxiv.org/pdf/1501.06164.pdf}. 

\bibitem[K10]{K9} N. Katzourakis,  \emph{Absolutely minimising generalised solutions to the equations of vectorial Calculus of Variations in $L^\infty$}, Calculus of Variations and PDE, in press.

\bibitem[K11]{K10} N. Katzourakis,  \emph{Equivalence between weak and $\mD$-solutions for symmetric hyperbolic first order PDE systems}, ArXiv preprint, \url{http://arxiv.org/pdf/1507.03042.pdf}. 

\bibitem[K12]{K11} N. Katzourakis,  \emph{Mollification of $\mD$-solutions to fully nonlinear PDE systems}, ArXiv preprint, \url{http://arxiv.org/pdf/1508.05519.pdf}. 

\bibitem[KP]{KP} N. Katzourakis,  T. Pryer, \emph{On the numerical approximation of $\infty$-Harmonic mappings}, Nonlinear Differential Equations and Applications  23 (6), 1 - 23 (2016).

\bibitem[KP2]{KP2} N. Katzourakis,  T. Pryer, \emph{Second order $L^\infty$ variational problems and the $\infty$-Polylaplacian}, ArXiv preprint, \url{https://arxiv.org/pdf/1605.07880.pdf}. 

\bibitem[KM]{KM} N. Katzourakis, J. Manfredi, \emph{Remarks on the Validity of the Maximum Principle for the $\infty$-Laplacian}, Le Matematiche, Vol. LXXI (2016) Ð Fasc. I, 63 - 74, DOI: 10.4418/2016.7 1.1.5.


\bibitem[KM2]{KM2} N. Katzourakis, R. Moser, \emph{Existence, uniqueness and structure of 2nd order absolute minimisers}, ArXiv preprint, \url{https://arxiv.org/pdf/1701.03348.pdf}. 

\bibitem[KR]{KR} J. Kristensen, F. Rindler, \emph{Characterization of generalized gradient Young measures generated by sequences in $W^{1,1}$ and $BV$}, Arch. Rational Mech. Anal. 197, 539 - 598 (2010) and \emph{erratum} Arch. Rational Mech. Anal. 203, 693 - 700 (2012).

\bibitem[P]{P} P. Pedregal, \emph{Parametrized Measures and Variational Principles}, Birkh\"auser, 1997.

\bibitem[S]{S} O. Savin, \emph{$C^1$ Regularity for Infinity Harmonic Functions in Two Dimensions}, Arch. Rational Mech. Anal. 176, 351 - 361, (2005).

\bibitem[SS]{SS} S. Sheffield, C.K. Smart, \emph{Vector Valued Optimal Lipschitz Extensions}, Comm. Pure Appl. Math. 65, 128 -154 (2012).

\bibitem[U]{U} K. Uhlenbeck, \emph{Regularity for a class of nonlinear elliptic systems}, Acta Mathematica 138, 219 - 240 (1977).

\bibitem[V]{V} M. Valadier, \emph{Young measures}, in ``Methods of nonconvex analysis", Lecture Notes in Mathematics 1446, 152-188 (1990).

\end{thebibliography}

\end{document}